\documentclass{procDDs}
\usepackage{xcolor,soul}
\usepackage{ulem,comment,hyperref}                                                  

\title{\vspace{-5.5mm}Computations of general Heun functions from their integral series representations\vspace{-2mm}}

\author{\textbf{T. Birkandan}}
{Department of Physics, Istanbul Technical University, 34469 Istanbul, Turkey}                 
{birkandant@itu.edu.tr}                                   
\coauthor{P.-L. Giscard}                       
{Univ. Littoral C\^ote d’Opale, UR 2597, LMPA, Laboratoire de Math\'ematiques Pures et Appliqu\'ees
Joseph Liouville, F-62100 Calais, France}                  
{giscard@univ-littoral.fr }                                               

\coauthor{A. Tamar}                       
{Independent researcher, India}                  
{adityatamar@gmail.com\vspace{-2mm}}                                               
                                                  
\begin {document}

\maketitle

\index{Birkandan, T.}                              
\index{Giscard, P.-L.}                              
\index{Tamar, A.}                             %

\begin{abstract}
We present a numerical implementation of the recently developed unconditionally convergent representation of general Heun functions as integral series. We produce two codes in Python available for download, one of which is especially aimed at reproducing the output of \textsc{Mathematica}'s HeunG function. 
We show that the present code compares favorably with Mathematica's HeunG and with an Octave/Matlab code of Motygin, in particular when the Heun function is to be evaluated at a large number of points if less accuracy is sufficient. We suggest further improvements concerning the accuracy and discuss the issue of singularities.
\end{abstract}

\section{Introduction}

Heun-type equations have found a large number of applications in physics, especially since the early 2000s \cite{Hort1, Hort2,Ishk4}.  
In a recent publication \cite{giscard2020}, a novel formulation of all functions of the Heun class was given as unconditionally convergent integral series.  While seemingly difficult to implement numerically owing to the large number of iterated integrals the method requires, we expose here significant simplifications allowing a fast and reliable mean to numerically evaluate the required integral series. 
The approach exposed here as well the underlying mathematical results \cite{giscard2020} are  valid for all functions of the Heun class: general, confluent, bi-confluent, doubly confluent and tri-confluent. For the sake of conciseness we present only the results pertaining to the general Heun function, namely,
 \begin{align}
&\frac{d^2H_G(z)}{dz^2} + \bigg[ \frac{\gamma}{z} + \frac{\delta}{z - 1} + \frac{\epsilon}{z - t} \bigg]\frac{dH_G(z)}{dz} 
\label{HeunEG}\\
&\hspace{25mm}+ \frac{\alpha\beta z - q}{z(z-1)(z-t)}H_G(z) = 0,\nonumber
\end{align}
where $t$ is the location of the regular singular point other than $\{0,\, 1,\, \infty\}$ \cite{SpecFunc,Ronveaux1995}. By following a confluence procedure on its singularities \cite{SpecFunc} other equations of Heun class are obtained from (\ref{HeunEG}).

\subsection{Mathematical background}
The approach developed in \cite{giscard2020} relies on the following observation: we can design a $2\times 2$ matrix $\mathsf{M}(z)$ depending on a complex variable $z$ such that the unique solution $\mathsf{U}(z,z_0)$ to the equation
\begin{equation}\label{DiffSys}
\frac{d}{dz}\mathsf{U}(z,z_0) = \mathsf{M}(z)\mathsf{U}(z,z_0),
\end{equation}
with initial condition $\mathsf{U}(z_0,z_0) = \mathsf{Id}$ the $2\times 2$ identity matrix, comprises a Heun function. More precisely, setting 
\begin{equation}\label{Mmatrix}
\mathsf{M}(z) = \begin{pmatrix}1&1\\
 B_1(z)+B_2(z)-1&B_1(z)-1\end{pmatrix},
\end{equation} 
with 
\begin{align*}
B_1(z) &:=  -\frac{\gamma}{z} - \frac{\delta}{z - 1} - \frac{\epsilon}{z - t},\\
B_2(z) &:= -\frac{\alpha\beta z - q}{z(z-1)(z-t)},
\end{align*}
and
\begin{equation}\label{IniVec}
\psi(z_0) :=\begin{pmatrix}H_0\\H_0'-H_0\end{pmatrix},
\end{equation}
then  $\psi(z) :=\mathsf{U}(z,z_0).\psi(z_0) $ is the unique general Heun function $H_G(z)$ solution of the Cauchy problem given by equation (\ref{HeunEG}) with boundary conditions $H_G(z_0)=H_0$ and $H_G'(z_0) = H_0'$. In itself this reformulation of the problem gains us nothing because matrix $\mathsf{M}(z)$ does not commute with itself at different $z$ values: $\mathsf{M}(z).\mathsf{M}(z')\neq \mathsf{M}(z')\mathsf{M}(z)$ when $z\neq z'$. Because of this $\mathsf{U}(z,z_0)$, called the evolution operator of system (\ref{DiffSys}), cannot easily be expressed in terms of $\mathsf{M}(z)$. Instead, $\mathsf{U}(z,z_0)$ is the ordered exponential of $\mathsf{M}(z)$, a generalisation of the notion of ordinary matrix exponential. 

Evaluating ordered exponentials analytically is an old and difficult problem. Until recently it was only possible to use either one of two approaches: Floquet-based methods which only provide approximate solutions and require $\mathsf{M}(z)$ to be periodic as a function of $z$ \cite{Floquet1883}; or the Magnus expansion \cite{Magnus1954}, a continuous analog of the Baker-Campbell-Hausdorff formula which involves an infinite series of increasingly intricate expressions and suffers from incurable divergence issues \cite{Blanes2009}. In 2015, \cite{Giscard2015} gave a graph-theoretic formulation of ordered exponentials allowing their exact expression in terms of continued fractions of finite depth and breadth taken with respect to a convolution like product, known as the Volterra composition.  

While it is not necessary to review the graph theoretic arguments underlying this result, it is most important for our purpose to recall the theory of Volterra composition. To this end, consider the space $D$ of distributions of the form
$$
f(z',z) = \tilde{f}(z',z)\Theta(z'-z)+\sum_{i=1}^\infty \tilde{f}_i(z',z)\delta^{(i)}(z'-z),
$$
where $\delta^{(i)}$ is the $i$th derivative of the Dirac delta distribution $\delta^{(0)}(z'-z)\equiv \delta(z'-z)$, $\Theta(z'-z)$ is the Heaviside step function with the convention $\Theta(0)=1$ and all $\tilde{f}_i(z',z)$ are smooth functions in both variables over an interval of interest $I$. We define the following product, for $f,l\in D$,
$$
(f\ast l)(z',z)=\int_{-\infty}^\infty f(z',\zeta)l(\zeta,z) d\zeta.
$$
This product makes $D$ into a non-commutative algebra \cite{GiscardPozzainv, GiscardPozzatri} with unit element $1_\ast\equiv \delta(z'-z)$. Now consider the ensemble $\text{Sm}_{\Theta}\subsetneq D$ of distributions of $D$ for which all smooth coefficients $\tilde{f}_i(z',z)=0$, i.e. only the Heaviside part remains. On $\text{Sm}_{\Theta}$ the $\ast$-product simplifies to 
\begin{align}
\hspace{-1mm}(f\ast l)(z',z)&\hspace{-.5mm}=\!\!\int_{-\infty}^\infty\!\!\!\! \tilde{f}(z',\zeta)\tilde{l}(\zeta,z) \Theta(z'-\zeta)\Theta(\zeta,z)\, d\zeta,\nonumber\\
&\hspace{-0.5mm}=\int_{z}^{z'}\!\! \tilde{f}(z',\zeta)\tilde{l}(\zeta,z) d\zeta\,\Theta(z'-z),\label{SmThetaProduct}
\end{align}
which is the convolution-like product introduced by Volterra in his studies of integral equations, now known as the Volterra composition \cite{Volterra1924}. Not only do distributions of $\text{Sm}_{\Theta}$ have $\ast$-inverses for all $z',z\in\mathbb{C}$ except at a countable number of isolated points \cite{GiscardPozzatri} but, most importantly for our purposes here, they have $\ast$-resolvents. Such resolvents, denoted $(1_\ast - f)^{\ast -1}$ are given by the Neumann series, 
 \begin{align*}
 \big(1_\ast - f\big)^{\ast -1}(z',z) &= \sum_{n=0}^\infty f^{*n}(z',z), 
 \end{align*}
 where $f^{\ast 0}=1_\ast$ and $f^{\ast i}=f\ast f^{\ast (i-1)}$. Writing the $\ast$-products explictly yields the integral series
 \begin{align*}
 &\big(1_\ast - f\big)^{\ast -1}(z',z)= \delta(z'-z)+\tilde{f}(z',z)\Theta(z'-z)\\ 
 &\hspace{1mm}+ \int_{z}^{z'} \tilde{f}(z',\zeta)\tilde{f}(\zeta,z)d\zeta\Theta(z'-z)\\
&\hspace{1mm}+\hdots
 \end{align*}
For all $f\in \text{Sm}_\Theta$ this series is unconditionally convergent everywhere on $z',z\in I^2$ except at the singularities of $f$ \cite{Giscard2015}. The $\ast$-resolvent of $f$ solves  
$$
f \ast \big(1_\ast - f\big)^{\ast -1} = \big(1_\ast - f\big)^{\ast -1}-1_\ast,
$$
which means that it solves the linear Volterra integral equation of the second kind with kernel $f(z',z) = \tilde{f}(z',z)\Theta(z'-z)$. The method of path-sum gives any entry of any ordered exponential exactly as continued fractions of $\ast$-resolvents \cite{Giscard2015}, each of which is representable by its unconditionally convergent Neumann series. 

For the matrix $\mathsf{M}$ of Eq.~(\ref{Mmatrix}) and initial vector of (\ref{IniVec}), the path-sum continued fraction yields
\begin{align}
  H_G(z)&=H_0+H_0\int_{z_0}^z \!G_{1}(\zeta,z_0)d\zeta\,\,+\label{HeunRes}\\
&\hspace{-7mm}(H'_0-H_0)\left(\!e^{z-z_0}-1+\!\int_{z_0}^z\!(e^{z-\zeta}-1)G_{2}(\zeta,z_0)d\zeta\right),\nonumber
  \end{align}
where the $G_i$ are defined from $\ast$-resolvents, 
$G_{i}=(1_\ast - K_i)^{\ast -1} -1_\ast = \sum_{n=1}^\infty K_{i}^{\ast n}$ with kernels  
   \begin{align*}
K_{1}(z,z_0)&=1+ \frac{e^{-z}}{z^{\gamma }
   (z-1)^{\delta } (t-z)^{\epsilon }}\times\\
   &\int_{z_0}^z \zeta_1^{\gamma } (\zeta_1-1)^{\delta }(t-\zeta_1)^{\epsilon }e^{\zeta_1} X(\zeta_1)d\zeta_1,\\
K_{2}(z,z_0)&=X(z) e^{z-z_0}-\frac{q-\alpha  \beta  z}{(z-1) z
   (z-t)}.
 \end{align*} 
 where 
 $$
 X(z):=\frac{q-\alpha  \beta  z}{(z-1) z
   (z-t)}-\frac{\epsilon
   }{t-z}-\frac{\gamma
   }{z}-\frac{\delta }{z-1}-1.
 $$
 See \cite{giscard2020} for the proof of this result. 
 \section{Numerical implementation}
 \subsection{Integral series on discrete times}
 Let $I = ]a,b[$ be the interval of interest over which the Heun function is sought, $z_0\in\bar{I}$ (where $\bar{I}$ designates the closure of $I$) be the point at which the initial values are known 
 and 
 let $\{z_i\in I\}_{0\leq i \leq N-1}$ be the discrete values at which the Heun function is to be numerically evaluated. For simplicity, suppose that the distance $|z_{i+1}-z_{i}|=\Delta z$ is the same for all $0\leq i\leq N-2$. This assumption is not necessary but will alleviate the notation. 
 Now for a smooth function $\tilde{f}(z',z)$ over $I^2$, we define a matrix $\mathsf{F}$ with entries 
 $$
 \mathsf{F}_{i,j} := f(z_i,z_j) = \tilde{f}(z_i,z_j)\Theta(z_i-z_j) .
 $$
 Note that, by construction, $\mathsf{F}$ is lower \textit{triangular} owing to the Heaviside step function.
 Of major importance is the observation that once $I$ is discretized, the Volterra composition of Eq.~(\ref{SmThetaProduct}) turns into an ordinary matrix product 
 \begin{align*}
 (f\ast l)(z_i,z_j) = &\int_{z_j}^{z_i} \tilde{f}(z_i,\zeta)\tilde{l}(\zeta,z_j) d\zeta~ \Theta(z'-z)\\
 &\hspace{10mm}\big\downarrow\\
 &\hspace{-16mm}\sum_{z_j\leq z_k\leq z_i} \!\!\tilde{f}(z_i,z_k)\tilde{l}(z_k, z_j)\, dz = (\mathsf{F}.\mathsf{L})_{i,j}\, \Delta z.
 \end{align*}
 Here $\mathsf{L}_{k,j}:=l(z_k,z_j)\Theta(z_k-z_j)$. Rigorously,
 $$
 \lim_{\Delta z\to 0}(\mathsf{F}.\mathsf{L})_{i,j}\, \Delta z = \int_{z_j}^{z_i} f(z_i,\zeta)l(\zeta,z_j)dz,
 $$
 This matricial presentation is useful for direct integration as well since for example
 \begin{equation}
 \int_{z_i}^{z_j} f(\zeta, z_i) dz\,\Theta(z_j-z_i)=\lim_{dz\to 0}  (\mathsf{H}.\mathsf{F})_{i,j}\,\Delta z,\label{intNum}
 \end{equation}
 where $\mathsf{H}_{i,j}:= 1$ if $i\geq j$ and 0 otherwise. 
 Most importantly, this line of results extends to $\ast$-resolvents, 
 $$
 \big(1_\ast - f\big)^{\ast -1}(z_i,z_j) =\lim_{\Delta z\to 0} \frac{1}{\Delta z}\big(\mathsf{Id}-\Delta z\,\mathsf{F}\big)^{-1}_{i,j},
 $$
 so that, in practical numerical computations with $\Delta z\ll 1$, we may  use
  $$
 (1_\ast - f)^{\ast -1}-1_\ast \simeq\big(\mathsf{Id}-\Delta z\,\mathsf{F}\big)^{-1}_{i,j}-\mathsf{Id}_{i,j}/\Delta z. 
 $$
 We improve on the above by noting that these results correspond to using the rectangular rule of integration as revealed by Eq.~(\ref{intNum}). Using the trapezoidal rule instead leads to much more accurate results.\footnote{In principle, arbitrary quadrature rules are possible but further work is needed, as determining the ensuing matrix forms for $\ast$-resolvents is not trivial.}
 To follow this rule, the usual matrix product $\mathsf{F}.\mathsf{L}\, \Delta z$ representing $f\ast l$ must be replaced by 
 $$
 \frac{\Delta z}{2}(\mathsf{F}-\mathsf{dF}).\mathsf{L} + \frac{\Delta z}{2}\mathsf{F}.(\mathsf{L}-\mathsf{dL}), $$
 where e.g. $\mathsf{dF}$ is the diagonal of $\mathsf{F}$. Now the $\ast$-resolvent of $f$ becomes
 $$
 (1_\ast - f)^{\ast -1}(z_i,z_j) \simeq \frac{1}{\Delta z}\left(\mathsf{Id} - \Delta z\mathsf{F} + \frac{\Delta z}{2}\mathsf{dF}\right)^{-1}_{i,j}.
 $$
 Taken together these results completely bypass the iterated integral series rather than evaluating directly what it converges to. The entire complexity inherent in the integral series representation of the general Heun function is reduced to multiplying and inverting \textit{triangular}, well-conditioned matrices. For the latter assertion, consider the diagonal elements of $\mathsf{Id}-\Delta z \mathsf{F}+(\Delta z/2)\mathsf{dF}$: these can be made to be different from  0 by tuning $\Delta z$.  
 
 \subsection{Numerical implementation}
 Concretely, a code implementing the path-sum formulation of the evolution operator $\mathsf{U}(z,z_0)$ and, from there, the general Heun function, must evaluate kernels $K_1$ and $K_2$, construct the corresponding matrices, compute their ordinary matricial resolvents, add and multiply them as required. 

 Let $\mathsf{K}_{1,2}$ be the triangular matrix representing kernel $K_1$ or $K_2$, i.e. $(\mathsf{K}_{1,2})_{i,j}=K_{1,2}(z_i,z_j)\Theta(z_i-z_j)$. Then, construct the matrix
 \begin{align*}
&\mathsf{R} := \\
&H_0\left(\mathsf{H} +\frac{\Delta z}{2} \mathsf{g}_1.(\mathsf{H}-\mathsf{dH})+\frac{\Delta z}{2}(\mathsf{g}_1-\mathsf{dg}_1).\mathsf{H}\right)+\\& 
(H_0'-H_0)\left( \mathsf{E} + \frac{\Delta z}{2} (\mathsf{E}-\mathsf{dE}).\mathsf{G}_2\frac{\Delta z}{2}+\mathsf{E}.(\mathsf{G}_2-\mathsf{dG}_2)\right),
 \end{align*}
 where $\mathsf{g}_1:=\mathsf{Id} +\Delta z \mathsf{G}_1$, $\mathsf{E}_{i,j} := (e^{z_i-z_j}-1)\Theta(z_i-z_j)$ and $\mathsf{G}_{1,2} := \Delta z^{-1}(\mathsf{Id} - \Delta z\mathsf{K}_{1,2} + (\Delta z/2)\mathsf{dK}_{1,2})^{-1}-\mathsf{Id}/\Delta z$. Comparing the above with Eq.~(\ref{HeunRes}) indicates that
 $$
\lim_{\Delta z\to 0} \mathsf{R}_{i,j} = \mathsf{U}(z_i,z_j).\psi(z_0),
 $$
 with an error scaling as $O(\Delta z^2)$, so that in particular
 $$
 \lim_{\Delta z\to 0} \mathsf{R}_{i,0} = H_G(z_i).
 $$
 This shows that only the first column of $\mathsf{R}$ is useful. This observation is profitably exploited numerically, as it avoids the need to even compute matrix inverses $(\mathsf{Id} - \Delta z\mathsf{K}_{1,2} + (\Delta z/2)\mathsf{dK}_{1,2})^{-1}$, rather asking for the solution of the triangular system $(\mathsf{Id} - \Delta z\mathsf{K}_{1,2} + (\Delta z/2)\mathsf{dK}_{1,2}).\vec{x}=\vec{v}$,
 where $v=(1,0,0...)^{\text{T}}$, which is faster and requires less memory. From there further matrix multiplications can be implemented vectorially on $\vec{x}$.
 
 The last remaining difficulty lies in constructing matrix $\mathsf{K}_1$ as $K_1(z_i,z_j)$ involves an integral
 $$
 \mathfrak{I}(z_i,z_j):=\int_{z_j}^{z_i} \zeta_1^{\gamma } (\zeta_1-1)^{\delta }(t-\zeta_1)^{\epsilon }e^{\zeta_1} X(\zeta_1)d\zeta_1.
 $$
 Consequently, computing all entries $(\mathsf{K}_1)_{i,j}$ naively requires a quadratically growing number of integrals $\mathfrak{I}(z_i,z_j)$ to be evaluated. This difficulty disappears thanks to the integral linearity: first 
 we compute only integrals $\mathfrak{I}(z_{j+1},z_j)$ over small intervals $[z_i,z_{i+1}]$, which is done e.g. with the trapezoidal rule. Then the first column of $\mathsf{K}_1$ is $$
 (\mathsf{K}_1)_{i,0} = 1+ \frac{e^{-z_i}}{z_i^{\gamma }(z_i-1)^{\delta } (t-z_i)^{\epsilon }} \mathfrak{I}(z_i,z_0),
 $$ 
 the required integrals of which we obtain from cumulative sums over the integrals on small intervals since $\mathfrak{I}(z_i,z_0)=\sum_{j=0}^{i-1}\mathfrak{I}(z_{j+1},z_j)$. The integrals $\mathfrak{I}(z_i,z_j)$ with $j>0$ which are required to construct the subsequent columns of $\mathsf{K}_1$ are similarly given  by $\mathfrak{I}(z_i,z_j)=\mathfrak{I}(z_i,z_0)-\mathfrak{I}(z_j,z_0)$, which can be done at once for all $j\geq 0$ with a meshgrid function.   
 
 At this point the numerical code relies solely on well-conditioned triangular matrices, need not inverting any of them explicitly as it solves two matrix-vector linear equations instead, computes no integral explicitly but relies on the trapezoidal rule throughout, is guaranteed to be everywhere convergent except at singularities and is exact under the limit $\Delta z\to 0$, with an expected error scaling of $O(\Delta z^2)$.   
 
 \subsection{Optional improvement}
 
 Seeking more flexibility in the calculations, the code allows one to subdivide the interval of interest $I$ into $N_1$ smaller intervals. Over each of the subintervals, $H_G(z)$ is evaluated from its path-sum formulation as presented in the preceding sections with $N_2$ points. The values of $H_G(z)$ and $H'_G(z)$ at the border of an interval is used as Cauchy initial value for the computations over the next interval. Overall, a total of $N=N_1N_2$ points are computed. Tuning $N_1$ and $N_2$ independently  allows the user to trade accuracy for large $N$ values and vice-versa. For example, picking $N_1\gg 1$ while keeping $N_2$ moderate allows for much faster evaluation than pure path-sum should the user require a very large number $N\gg 1 $ of values, or to evaluate $H_G$ over a very large interval. At the opposite, more accurate results will be obtained by making $N_2$ large while $N_1$ can be as low as 1.

 \section{Performance}
 \subsection{Comparisons with Mathematica's $HeunG$ and Octave's HeunL0}\label{Sec:Compar}
  \textsc{Mathematica}'s HeunG function is a new feature in version 12.1 of the software. It produces the regular solution of the general Heun equation with boundary conditions $H_G(0) = 1$ and $H_G'(0)=q/(\gamma\,t)$. This quantity is also produced by function HeunL0 from Motygin's general purpose Heun package for \textsc{Octave}/\textsc{MATLAB} \cite{Motygin2015}. We propose to compare the performances of these two codes with that of the integral-series based Python code available on GitHub \cite{PythonCode}.
   \begin{table}
    \begin{center}
  \vspace{-2mm}
  \begin{small}
      
\begin{tabular}{| c| c| c |c|}
\hline\hline
 Code  & Accuracy& Points & Time (sec.)\\
 \hline
 HeunG &$10^{-16}$& 1000 & 2.22 \\
 HeunL0 &$10^{-13}$& 1000&6.05\\
 HeunL0 &$10^{-6}$& 1000&2.19\\
 This code  &$10^{-6}$& 1000 & 0.0096 \\
 \hline
 HeunG &$10^{-16}$& 10000 & 22.4 \\
 HeunL0&$10^{-13}$& 10000& 62.26\\
 HeunL0 &$10^{-6}$& 10000&21.80\\
 This code  &$10^{-6}$& 10000 & 0.090 \\
 \hline
HeunG &$10^{-16}$& 50000 & 110.8 \\
HeunL0&$10^{-13}$& 50000& 297.99\\
HeunL0 &$10^{-6}$& 50000&109.36\\
 This code  &$10^{-6}$& 50000 & 0.43 \\
 \hline
 HeunG &$10^{-16}$& 100000 & 231 \\
 HeunL0&$10^{-13}$& 100000&649.15\\
 HeunL0 &$10^{-6}$& 100000&219.19\\
 This code  &$10^{-6}$& 100000 & 0.90 \\
 \hline
 HeunG &$10^{-16}$& 200000 & 464 \\
 HeunL0&$10^{-13}$& 200000&1274.6\\
 HeunL0 &$10^{-6}$& 200000&434.76\\
 This code  &$10^{-6}$& 200000 & 1.87 \\
 \hline
 \hline
 \end{tabular}
 \caption{\label{Tab:Compar} Comparisons of performances between codes computing a general Heun function, for parameters see \S\ref{Sec:Compar}.}
 \vspace{-2mm}
 \end{small}
\end{center}
\end{table}
  
  We chose the following parameters for the Heun function : $t = 9/2$, called `$a$' in \textsc{Mathematica}'s HeunG, $q=-1$, $\alpha=1$, $\beta=-3/2$, $\gamma=-14/100$, $\delta = 432/100$ and $\epsilon= 1.0 + \alpha + \beta - \gamma - \delta $. We want to represent $H_G(z)$ for $z\in I=[-2.2,0.8]$. Observe how this interval comprises a singularity at $z = 0$. Problems with this singularity are avoided upon specifying ordinary conditions on the left $z<0$ and right $z>0$ of it at $|z_0|\ll 1$. Values for $H_G(z_0)$ and $H_G'(z_0)$ are known  within machine precision from the first three orders of the general series expansion of $H_G(z)$ near $z=0$. For the more general issue of crossing over singularities, see the discussion in \S\ref{Singularities} below.
  \begin{figure*}[!h]
\centering
 \vspace{-5mm}
\includegraphics[width=.49\textwidth]{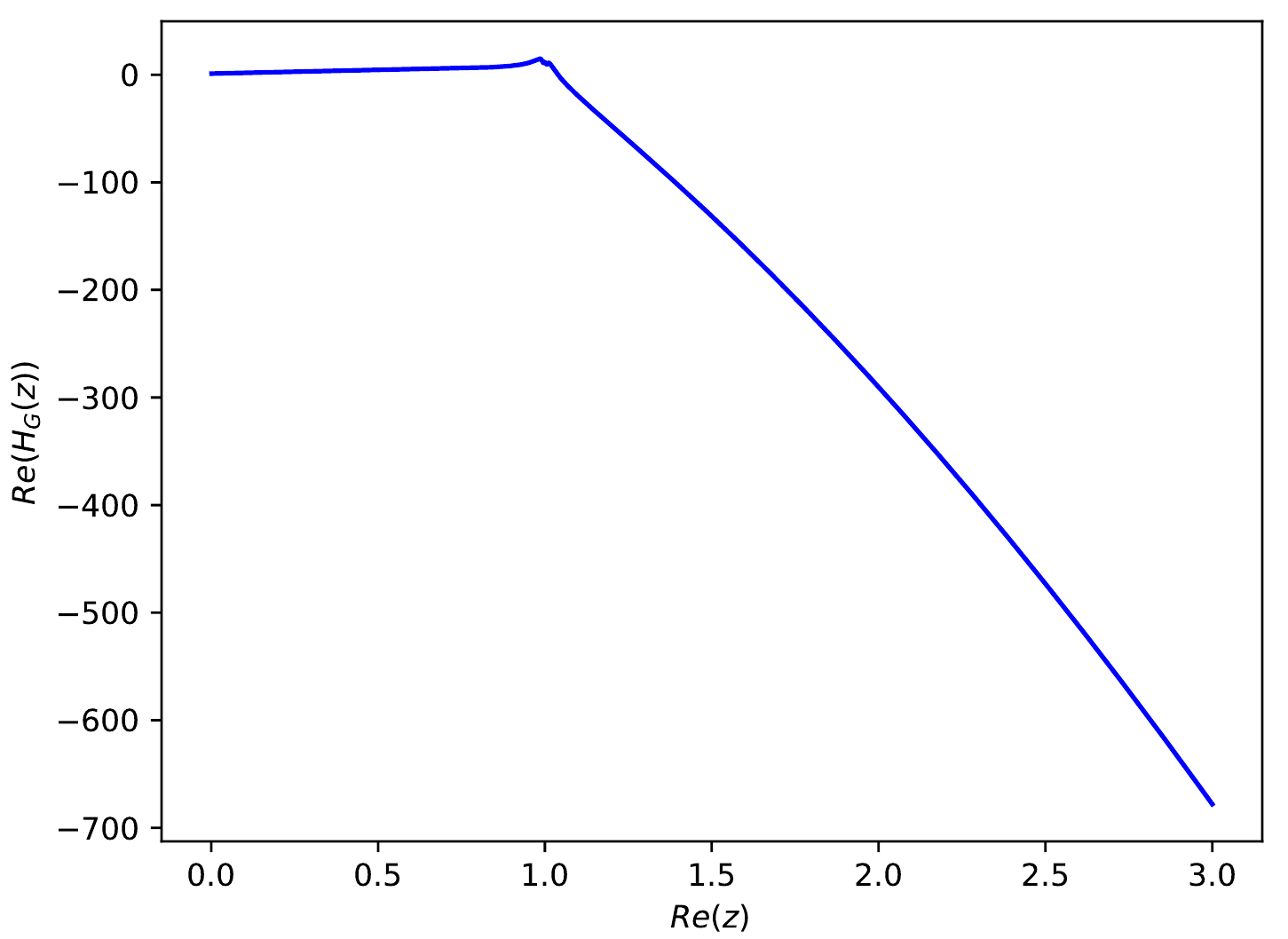}
\includegraphics[width=.485\textwidth]{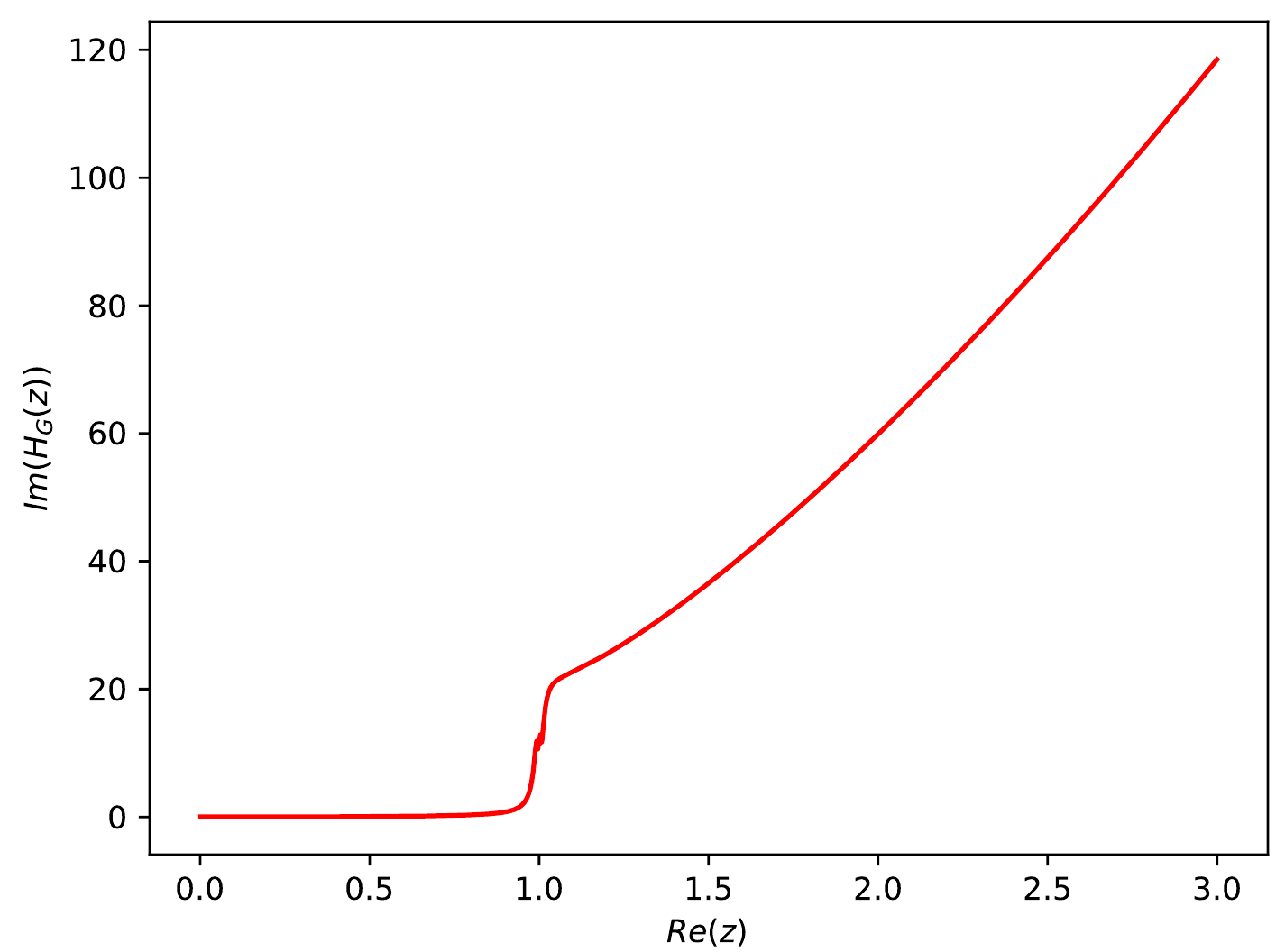}\\
\hspace{2.2mm}\includegraphics[width=.475\textwidth]{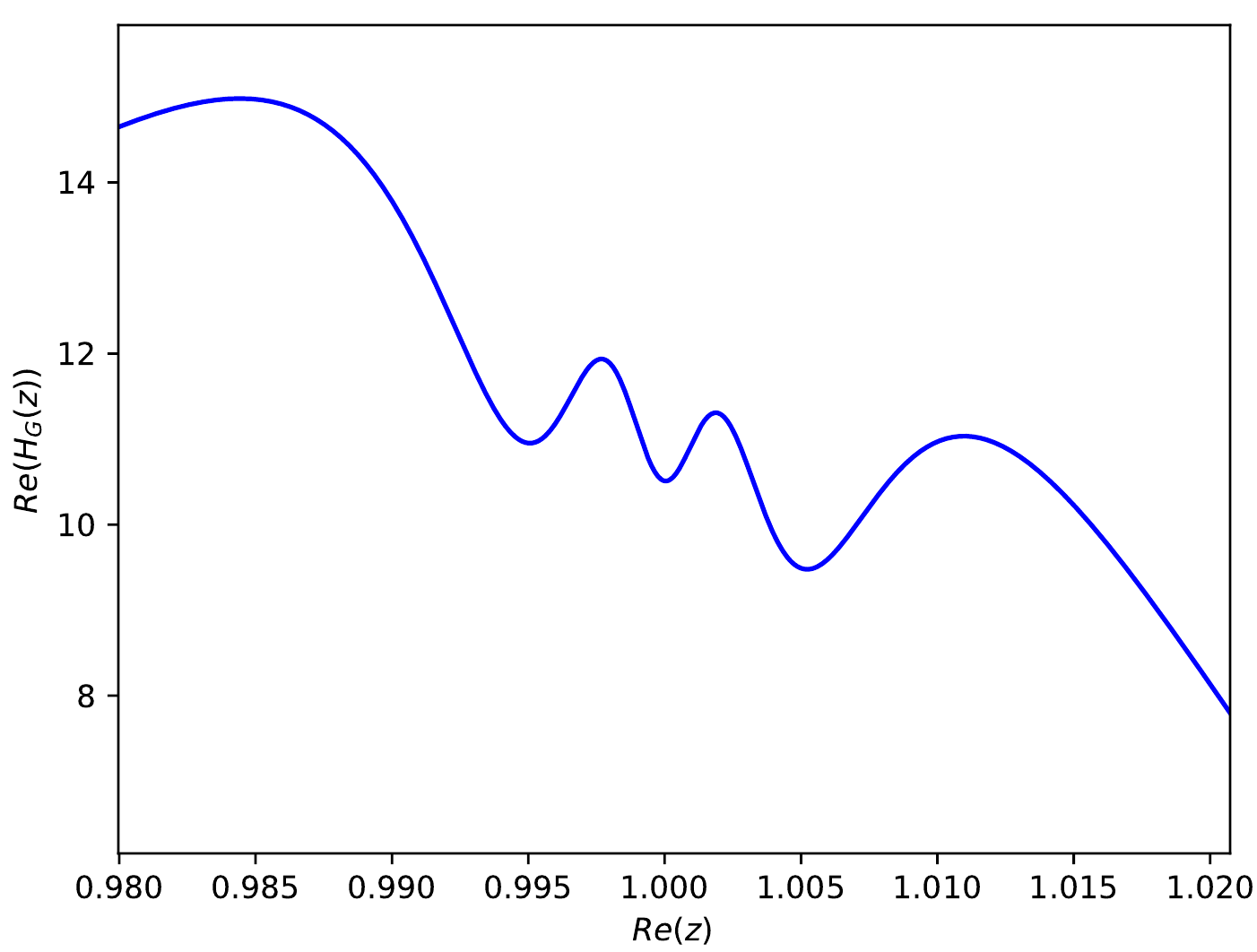}
\includegraphics[width=.485\textwidth]{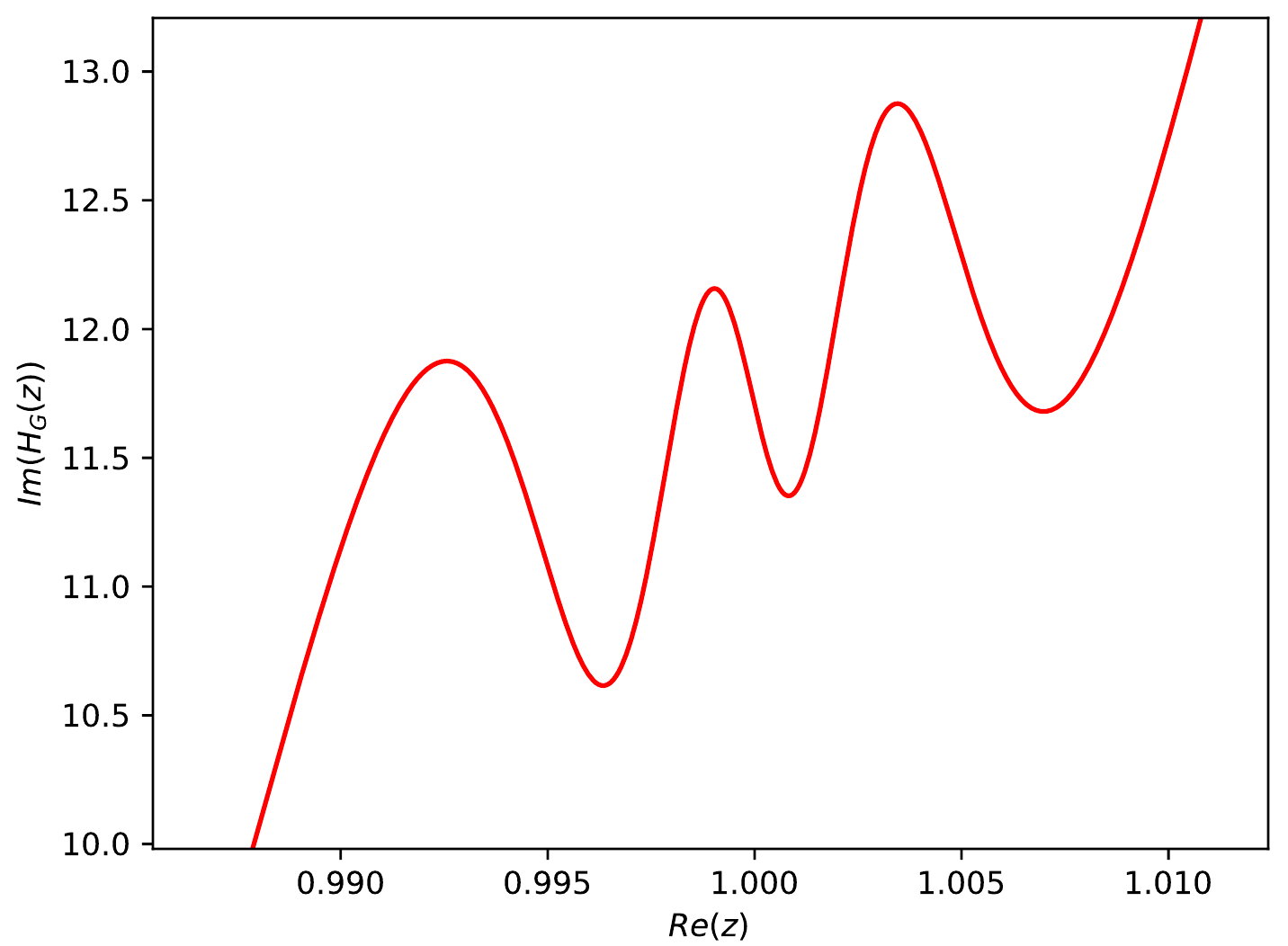}
\vspace{-2mm}
\caption{
General Heun function for $t = 1+i\times 10^{-2}$, $q=-1$, $\alpha=1$, $\beta=-3/2$, $\gamma=-14/100$, $\delta = 432/100$ and $\epsilon= 1.0 + \alpha + \beta - \gamma - \delta $ along the path $z=z_r+5i\times 10^{-3}$ with $z_r\in[0,3]$. Top row: real part $\Re(H_G(z))$ (left fig.) and imaginary part $\Im(H_G(z))$ (right fig.). Bottom row: close-ups of $\Re(H_G(z))$ (left fig.) and $\Im(H_G(z))$ (right fig.) near the singularities at $z=1$ and $z=t$.\label{HGSing}}
\vspace{-3mm}
\end{figure*}
All codes were tasked determining $H_G(z)$ for $z\in[-2.2,0.8]$ by producing a table of values $H_G(z_j)$ with $z_j = -2.2+ j\times 3/N$, where $N$ is the number of points. The Python code (`num$\_$heunG.py' available for download) used parameters $N_2=100$ and increasing $N_1$ so as to obtain the desired number $N=N_1N_2$ of points. This provides an experimentally observed accuracy of up to $10^{-6}$ for all points. This means the comparisons effected with HeunG are unfair in terms of accuracy as it works with machine precision. For HeunL0, we present two versions, one reaching $10^{-13}$ accuracy, the other one modified so as to reach an accuracy of $10^{-6}$ for fair comparisons with both HeunG and the Python code. 
 The Python code is written using Python 3.8. We employed the `time' module in Python, command `Timing[...]' in \textsc{Mathematica} and `tic;...toc;' in \textsc{Octave} to measure the run time. In this analysis, we used a Dell laptop running Ubuntu 18.04, equipped with Intel Core i7-8665U CPU @ 1.90GHz $\times 8$ and 16 GB memory, using a single core and serial computations. See Table.~\ref{Tab:Compar} for the results.

 We now turn to a challenging numerical problem, namely that of producing a general Heun function along a path in the complex plane $z = z_R + 5i\times10^{-3} $, $z_R\in[0,3]$, passing extremely close to two singularities at $1$ and $t=1+i\times 10^{-2}$. In the vicinity of $z\sim 1$ the Heun function changes abruptly (see Fig.~\ref{HGSing}) and we expect computations to experience numerical instabilities. As a consequence, we observe e.g. that HeunL0 requires well over the first 1000 terms in its series representation to estimate $H_G(z)$ at any one point for which $z_R\gtrsim 1$, reaching up to order $1542$ in $z_R=3$. In contrast, before the singularities $z_R\in[0,1[$, HeunL0 requires orders less than 50. HeunG becomes cripplingly slow when $z_R\sim 1$ to the point that we simply could not evaluate $H_G(z)$ in this region with \textsc{Mathematica}. The computation times for HeunL0 quickly becomes prohibitive for large numbers of points, already taking $563$ seconds for $10,000$ points so that we expect it to take over 7 hours for $\sim 500,000$ points. 
 The present code took $434$ seconds to produce $494,900$ points $(N_1=100, N_2=5000)$ with an accuracy of $\sim 10^{-6}$. It took $4.18$ seconds to produce the plots of Fig.~(\ref{HGSing}) on $49,400$ points $(N_1=100, N_2=500)$ with an accuracy $\sim 10^{-3}$.

 \subsection{General Heun functions with arbitrary boundary conditions}
The Python implementation of the integral series formulation of Heun function is natively well suited to evaluating the general Cauchy problem of determining the solution $H_G(z)$ of Eq.~(\ref{HeunEG}) with arbitrary Cauchy conditions at $z_0$ on an interval $I$ that does not cross over a singularity. The Python code (`num$\_$heun.py' available for download) accepts the parameters of the Heun function, the interval $I$ and the values of $N_1$ and $N_2$ as inputs.  

\subsection{Crossing singularities and exploring the complex plane}\label{Singularities}
The integral representation of Heun functions 
is valid and guaranteed to converge everywhere in the complex plane $z\in\mathbb{C}$ except at the singularities themselves. A difficulty thus presents itself when attempting to pass through a singularity with ordinary \emph{line} integrals, which fail to be defined there. Instead, we may integrate along any smooth curve $\gamma$ starting at $z_0$, ending at the desired $z_1$ and avoiding all singularities in between. Since the punctured plane is connected this will recover values for $H_G(z_1)$ and $H_G'(z_1)$ irrespectively of whether or not singularities lie on the line joining $z_0$ and $z_1$. These values are not mathematically anymore however due to the presence of branch cuts, which must be taken into account \cite{Motygin2015}. Keeping these caveats in mind, the numerical mathematics and the code presented here work the same replacing
$
\int_{z_0}^{z_1} d\zeta \mapsto \int_{\gamma\subset \mathbb{C}\backslash\mathcal{S}\atop \gamma:z_0\to z_1}d\gamma 
$
where $\mathcal{S}$ denotes the set of singularities of $H_G$. 
The concrete implementation of this extension to contours, while it does not seem to present fundamental difficulty, has yet to be implemented.   

%
\vspace{-4mm}

\section*{Acknowledgements}
\vspace{-2mm}
We thank the referee for constructive comments and suggestions. P.-L. G. is supported by the Agence Nationale de la Recherche young researcher grant No. ANR-19-CE40-0006.
\vspace{-4mm}

\bibliographystyle{unsrt}

\end{document}